\documentclass[11pt]{amsart}	

\usepackage[utf8]{inputenc}
\usepackage{graphicx}
\usepackage{subcaption}
\usepackage{amsmath}
\usepackage{amssymb}
\usepackage{amsthm}
\usepackage{ifthen}
\usepackage{color}

\def\polk#1{\setbox0=\hbox{#1}{\ooalign{\hidewidth\lower1.5ex\hbox{`}\hidewidth\crcr\unhbox0}}}

\newcommand{\osc}{\operatorname{osc}}

 \newcommand{\Rr}{\mathbb R}

\newcommand{\tr}{\operatorname{Tr}}

\newtheorem{Theorem}{Theorem}
\newtheorem{Definition}{Definition}
\newtheorem{Lemma}{Lemma}

\newtheorem{Proposition}{Proposition}
\newtheorem{assump}{}

\theoremstyle{definition}

\theoremstyle{remark}
\newtheorem{Remark}{Remark}

\usepackage{setspace}
\numberwithin{equation}{section}

\title[Fully nonlinear Hamilton-Jacobi equations]{Fully nonlinear Hamilton-Jacobi equations\\ of degenerate type}

\author[D. Jesus]{David Jesus}
\address{University of Coimbra, CMUC, Department of Mathematics, 3001-501 Coimbra, Portugal}{}
\email{djbj@mat.uc.pt}

\author[E. A. Pimentel]{Edgard A. Pimentel}
\address{University of Coimbra, CMUC, Department of Mathematics, 3001-501 Coimbra, Portugal and Pontifical Catholic University of Rio de Janeiro -- PUC-Rio, 22451-900 G\'avea, Rio de Janeiro-RJ, Brazil}{}
\email{edgard.pimentel@mat.uc.pt}

\author[J.M.~Urbano]{Jos\'{e} Miguel Urbano}
\address{King Abdullah University of Science and Technology (KAUST), Computer, Electrical and Mathematical Sciences and Engineering Division (CEMSE), Thuwal 23955-6900, Saudi Arabia and University of Coimbra, CMUC, Department of Mathematics, 3001-501 Coimbra, Portugal}{} 
\email{miguel.urbano@kaust.edu.sa} 

\date{\today}

\begin{document}

\begin{abstract}
We examine Hamilton-Jacobi equations driven by fully nonlinear degenerate elliptic operators in the presence of superlinear Hamiltonians. By exploring the Ishii-Jensen inequality, we prove that viscosity solutions are locally Lipschitz-continuous, with estimates depending on the structural conditions of the problem. We close the paper with an application of our findings to a two-phase free boundary problem.
\end{abstract}

\dedicatory{To the fond memory of Emmanuele DiBenedetto.}


\maketitle

\section{Introduction}

Among other seminal contributions, Emmanuele DiBenedetto revolutioni\-sed the analysis of singular and degenerate elliptic and parabolic equations. The development of intrinsic scaling methods \cite{DiBe93, DiBeUrbVes04, Urb08, DiBeGiaVes11}, in particular, had a huge impact in the field, with far-reaching applications that resonate to this day. 

In this paper, we study a fully nonlinear Hamilton-Jacobi equation of the form
\begin{align}\label{eq_main}
   F(D^2u)+H(Du,x)=f(x)\hspace{.2in}\mbox{in}\hspace{.2in}\Omega\subset\mathbb{R}^d,
\end{align}
where $F:S(d)\to\mathbb{R}$ is degenerate elliptic, the Hamiltonian $H=H(p,x)$ satisfies natural growth and continuity conditions, and $f\in L^\infty(\Omega)$ is Lipschitz continuous. In the superlinear setting, we prove that viscosity solutions to \eqref{eq_main} are locally Lipschitz-continuous. In addition, we examine a two-phase free boundary problem driven by the operator in \eqref{eq_main}. In this context, our findings include the existence of solutions and regularity estimates across the free boundary. The conditions we impose on the structure of the problem are fairly general and cover important models, such as Bellman and Isaacs equations. An example of Hamiltonian falling under our assumptions is
\[
	H(p,x):=a(x)\left(1+|p|^2\right)^\frac{m}{2}+V(x),
\]
provided $a,V:\Omega\to\mathbb{R}$ are Lipschitz-continuous and bounded from above and below, and $m>1$. 

Hamilton-Jacobi (HJ) equations of second-order often relate to stochastic optimal control problems \cite{Fleming-Soner_1991}. In this context, the value function of the optimization problem is a viscosity solution of the associated HJ equation. For developments at the intersection of viscosity solutions, Hamilton-Jacobi equations and optimal control, we refer the reader to \cite{Crandall-Lions_1983,Crandall-Evans-Lions_1984,Ishii_1989,Ishii-Lions_1990,Barles_1994, Bardi-Dolcetta_1997}. Recent results in regularity theory for degenerate problems, both in the variational and non-variational settings, can be found in \cite{Teixeira-Urbano_2014, Teixeira_2014, Araujo-Teixeira-Urbano_2017, Huaroto-Pimentel-Rampasso-Swiech_2020}.

The study of \eqref{eq_main} in the case $F\equiv \tr$ appears, for instance, in \cite{Lions_1985}. In that paper, the author proves the existence of classical solutions for the problem under Neumann boundary conditions and natural growth regimes on the Hamiltonian $H$. The role of Neumann (or oblique) boundary conditions relates to state-constrained optimal control problems, as they encode a reflection at the boundary; we refer the reader to \cite{Lions-Sznitman_1984}.

An in-depth account of state-constrained optimal control problems is the subject of \cite{Lasry-Lions_1989}, where the authors examine \eqref{eq_main} in the uniformly elliptic setting. In addition to establishing local Lipschitz continuity of the solutions, the authors connect the characterisation of boundary conditions with the growth regime of the Hamiltonian. Indeed, for sub-quadratic Hamiltonians, solutions blow up as they approach $\partial\Omega$, requiring a relaxed notion of boundary condition. However, in the (strictly) superquadratic case, solutions are globally H\"older continuous. See \cite{Barles-Lio_2004} for related developments.

The regularity theory available for \eqref{eq_main} in the degenerate elliptic case advanced substantially with the contributions in \cite{Dolcetta-Leoni-Porretta_2010}. Among the findings in that paper, we highlight the H\"older continuity of \emph{sub-solutions} in the presence of superquadratic Hamiltonians. The remarkable aspect of this result is in the one-sided requirement entailed by the sub-solution condition. To properly appreciate the minimality of such an assumption, we briefly recall the Krylov-Safonov theory. Concerning uniformly elliptic fully nonlinear operators, the latter implies H\"older continuity provided a \emph{two-sided control} is available. Indeed, let $\mathcal{P}^\pm$ be the Pucci extremal operators and $C>0$ be a constant; if $u\in C(\Omega)$ is a viscosity solution to
\[
	\mathcal{P}^-(D^2u)\leq C\hspace{.2in}\mbox{in}\hspace{.2in}\Omega
\]
and
\[
	\mathcal{P}^+(D^2u)\geq -C\hspace{.2in}\mbox{in}\hspace{.2in}\Omega,
\]
the Krylov-Safonov theory ensures the H\"older continuity of $u$. However, if one of the former inequalities fails to hold, the theory is no longer available; see, for instance, \cite{Caffarelli-Cabre_1995, Mooney_2019, Pimentel_2022}. We notice the results in \cite{Dolcetta-Leoni-Porretta_2010} also include global H\"older continuity for sub-solutions and Lipschitz regularity for the solutions of the homogeneous problem. We also refer the reader to the important contributions in \cite{Barles_1991} and \cite{Barles_2010}.

In \cite{Armstrong-Tran_2015}, the authors study \eqref{eq_main}, considering $F(M,x)=\tr (A(x)M)$, for a degenerate elliptic matrix-valued map $A:\Omega\to\mathbb{R}^{d^2}$. The findings in \cite{Armstrong-Tran_2015} advance the general theory of Hamilton-Jacobi equations, as they cover a general maximum principle, Lipschitz continuity for the solutions of the homogeneous equation with explicit estimates (in terms of the matrix $A$ and the structure of $H$), and state-constrained boundary conditions. We notice the Lipschitz continuity result in \cite{Armstrong-Tran_2015} relies on the maximum principle \cite[Theorem 3.2]{Crandall-Ishii-Lions_1992} and explores the connection of the trace operator and eigenvalues.

Our contribution is two-fold. By developing an intrinsically nonlinear argument, we prove that viscosity solutions to \eqref{eq_main} are Lipschitz continuous, with estimates. Then we examine a consequence of our regularity result to a two-phase free boundary problem and prove the existence of solutions, with estimates in H\"older spaces. Our main result reads as follows.

\begin{Theorem}[Improved regularity of solutions]\label{teo_main1}
Let $u\in C(\Omega)$ be a viscosity solution to \eqref{eq_main} where $F:S(d)\to\mathbb{R}$ is degenerate elliptic, Lipschitz-continuous and positively homogeneous of degree $1$, and $f\in L^\infty(\Omega)$ is Lipschitz-continuous. Suppose the Hamiltonian $H$ is superlinear and satisfies natural growth and continuity conditions, detailed in Section \ref{sec_immunization}. Then $u$ is locally Lipschitz-continuous in $\Omega$. Moreover, for every $\Omega'\Subset\Omega$, there exists $C>0$ such that
\[
	\left|u(x)-u(y)\right|\leq C\left( 1+\left\|u\right\|_{L^\infty(\Omega)}+\left\|f\right\|_{L^\infty(\Omega)}\right)|x-y|,
\]
for every $x,y\in\Omega'$. The constant $C>0$ depends only on $\Omega'$, and the data of the problem.
\end{Theorem}

The proof of Theorem \ref{teo_main1} relies on two building blocks. First, we examine the superquadratic case, i.e., $H(p,x)\sim C+C|p|^m$, with $m>2$. In this setting, we  prove that solutions to \eqref{eq_main} are $\gamma$-H\"older continuous for $\gamma:=(m-2)/(m-1)$. Here we follow closely the strategy put forward in \cite[Section 3]{Armstrong-Tran_2015}, adapting its techniques to the fully nonlinear setting. Then we refine the application of the Ishii-Jensen Lemma to produce Lipschitz regularity in the superlinear setting.

Once Theorem \ref{teo_main1} is available we turn to a free boundary problem driven by a particular instance of the operator in \eqref{eq_main}. For constants $\lambda_+,\lambda_-\in\mathbb{R}$, consider the problem
\begin{equation}\label{eq_touriga}
	{\rm Tr}(A(x)D^2u)+H(Du,x)=\lambda_+\chi_{\{u>0\}}+\lambda_-\chi_{\{u<0\}}\hspace{.1in}\mbox{in}\hspace{.1in}\Omega(u),
\end{equation}
where $\Omega(u):=\{u>0\}\cup\{u<0\}$. We note that \eqref{eq_touriga} holds \emph{only} in the region where the solutions do not vanish, and no PDE information is available in $\{u=0\}$. In many cases, viscosity solutions exist following Perron’s method, when a comparison principle is available, and one can build appropriate sub- and supersolutions. Meanwhile, the present setting introduces important difficulties. First, the dependence of \eqref{eq_touriga} with respect to the solution implies the lack of properness. As a consequence, one should not expect a comparison principle to hold at the level of the equation, which precludes the use of usual arguments. Also, and perhaps even more important, the growth regime of the Hamiltonian $H$ requires further compatibility conditions for the boundary data; see \cite{Dolcetta-Leoni-Porretta_2010}. 

We argue through a regularization of the right-hand side of \eqref{eq_touriga}, removing the dependence of the equation on zero-order terms. We combine our findings in regularity theory with former results on the existence of solutions for superquadratic Hamilton-Jacobi equations. Then a fixed-point argument ensures the existence of viscosity solutions to the Dirichlet problem associated with \eqref{eq_touriga}. For a similar approach in the context of free transmission problems, see \cite{Pimentel-Swiech_2022}.

The remainder of this paper is organised as follows. Section \ref{sec_immunization} gathers our primary assumptions and recalls a few preliminaries. The proof of Theorem \ref{teo_main1} is the subject of Section \ref{sec_lipschitz}. Finally, in Section \ref{cleese}, we prove the existence of viscosity solutions for the Dirichlet problem associated with \eqref{eq_touriga}.

\section{Preliminaries}\label{sec_immunization}

Here we detail the main assumptions used in the paper and collect a few preliminaries.

\subsection{Main assumptions}\label{subsec_firstshot}

We proceed with the conditions imposed on the second-order operator $F$. 

Because we rely on the monotonicity of $F$ in the space of symmetric matrices, we equip the latter with a partial order relation. For $M,N\in S(d)$, we say that $M\geq N$ if $M-N$ is positive semi-definite, i.e., for every $\xi\in\mathbb{R}^d$, we have 
\[
	\xi^T(M-N)\xi\geq 0.
\]
Our primary condition on $F$ concerns its degenerate ellipticity. 

\begin{Definition} 
We say $F:S(d)\to\mathbb{R}$ is degenerate elliptic if 
$$F(M)\leq F(N)$$
whenever $M, N\in S(d)$ are such that $M\geq N$.
\end{Definition} 

We continue with an assumption combining degenerate ellipticity and Lipschitz continuity for $F$.

\begin{assump}[Monotonicity and Lipschitz continuity]\label{assump_deg_ellip}
 The operator $F:S(d)\to \mathbb{R}$ is monotone non-increasing and Lipschitz continuous. That is, there exists a constant $C_F>0$ such that 
\begin{equation}\label{eq_aaronshore}
	F(M)-F(N)\leq C_F\left|\left(N-M\right)_+\right|,
\end{equation}
for every $M,N\in S(d)$.
 \end{assump}
 
 The condition in \eqref{eq_aaronshore} is equivalent to requiring $F$ to be degenerate elliptic and Lipschitz continuous, with constant $C_F>0$; see Lemma \ref{lem_murray}. The choice for \eqref{eq_aaronshore} stems from our argument since we compare $F(M)$ and $F(N)$ in terms of the eigenvalues of $M$ and $N$. We also require $F$ to be positively homogeneous of degree one.

\begin{assump}[Homogeneity of $F$]\label{assump_homog_F}
The operator $F$ is positively homogeneous of degree $1$. That is, for every $M\in S(d)$ and every $s\geq 0$, we have 
\[
	F(sM)=sF(M).
\]
\end{assump}

The typical example of an operator satisfying assumptions \ref{assump_deg_ellip} and \ref{assump_homog_F} is the Bellman operator. Indeed, let $\mathcal{A}$ be a measurable index set and consider a family of matrices $(A_\alpha)_{\alpha\in\mathcal{A}}$ such that
\[
	0\leq A_\alpha\leq (C_Fd^{-1})I,
\]
for every $\alpha\in\mathcal{A}$. Then the operator 
\[
	F(M):=\inf_{\alpha\in\mathcal{A}}\left(-\tr\left(A_\alpha M\right)\right)
\]
satisfies both \ref{assump_deg_ellip} and \ref{assump_homog_F}. 

Now, we turn to the Hamiltonian $H$ and detail the growth and continuity conditions under which we work.

\begin{assump}[Structural conditions]\label{assump_H_new}
We suppose there exist constants $m>1$ and $C_1, C_2, C_3>0$ such that $H=H(p,x)$ satisfies
\begin{equation}\label{yafo}
	-C_1+C_2 |p|^m \leq H(p,x)\leq C_3\left(1+|p|^m\right),
\end{equation}
for every $p\in\mathbb{R}^d$ and $x\in \Omega$. Also,
\begin{equation}\label{bneibrak}
	\left| H(p,x)-H(p,y)\right|\leq \left(C_3|p|^m+C_1\right)\left|x-y\right|
\end{equation}
for every $x,y\in \Omega$, and $p\in\mathbb{R}^d$. Finally, we require
\begin{equation}\label{beersheva}
	\left| H(p,x)-H(q,x)\right|\leq C_3\left(|p|+|q|+1\right)^{m-1}\left|p-q\right|,
\end{equation}
for every $p,q\in\mathbb{R}^d$, and every $x\in \Omega$. 
\end{assump}

The typical example of a Hamiltonian satisfying \ref{assump_H_new} is
\[
	H(p,x)=a(x)\left(1+|p|^2\right)^\frac{m}{2}+V(x),
\]
where $a,V:\Omega\to\mathbb{R}$ are Lipschitz continuous, with
\[
   0<C_\ast\leq a(x)\leq C^\ast\hspace{.2in}\mbox{and}\hspace{.2in }0\leq V(x)\leq C^\ast,
\]
for some fixed constants $0<C_\ast\leq C^\ast$.

\subsection{Preliminary material}\label{subsec_secondshot}

We start with the definition of viscosity solution.

\begin{Definition}[Viscosity solution]\label{israelidates}
We say $u\in {\rm USC}(\Omega)$ is a viscosity sub-solution to \eqref{eq_main} if, for every $x_0\in\Omega$ and every $\varphi\in C^2(\Omega)$ such that $u-\varphi$ has a local maximum at $x_0$, we have
\[
	F(D^2\varphi(x_0))+H(D\varphi(x_0),x_0)\leq f(x_0).
\]
Likewise, we say that $u\in {\rm LSC}(\Omega)$ is a viscosity supersolution to \eqref{eq_main} if,  for every $x_0\in\Omega$ and every $\varphi\in C^2(\Omega)$ such that $u-\varphi$ has a local minimum at $x_0$, we have
\[
	F(D^2\varphi(x_0))+H(D\varphi(x_0),x_0)\geq f(x_0).
\]
If $u\in C(\Omega)$ is simultaneously a viscosity sub-solution and a viscosity supersolution to \eqref{eq_main}, we say it is a viscosity solution to \eqref{eq_main}.
\end{Definition}

Next, we recall a maximum principle available for degenerate elliptic operators, see \cite[Theorem 3.2]{Crandall-Ishii-Lions_1992}.

\begin{Proposition}[Ishii-Jensen Lemma]\label{prop_joaquin}
Let $\Omega$ be a bounded domain and $G, \, H\in C\left(\Omega\times \Rr^d\times S(d)\right)$ be degenerate elliptic operators. Let $u_1$ be a viscosity sub-solution of $G\left(x, Du_1, D^2u_1\right)=0$ and $u_2$ be a viscosity supersolution of $H\left(x, Du_2, D^2u_2\right)=0$ in $\Omega$. Define $v:\Omega\times \Omega\to \Rr$ by
\[
    v\left(x,y\right):=u_1(x)-u_2(y)
\]
and  let $\varphi\in C^2(\Omega\times \Omega)$. Suppose that $\left(\overline{x}, \overline{y}\right)\in\Omega\times \Omega$ is a local maximum for $v-\varphi$. Then for every $\varepsilon>0$, there exist matrices $X_\varepsilon$ and $Y_\varepsilon$ in $S(d)$ such that
\[
    G\left(\overline{x}, D_x\varphi\left(\overline{x}, \overline{y}\right), X_\varepsilon\right)\leq 0 \leq H\left(\overline{y}, -D_y\varphi\left(\overline{x}, \overline{y}\right), Y_\varepsilon\right).
\]
Moreover, the matrix inequality 
\[
    -\left(\frac{1}{\varepsilon}+\left\|J\right\|\right)I\leq 
    \begin{pmatrix}
  X_\varepsilon & 0 \\
  0 & -Y_\varepsilon 
 \end{pmatrix}
 \leq J+\varepsilon J^2
\]
holds true, where $J:=D^2\varphi\left(\overline{x}, \overline{y}\right)$.
\end{Proposition}

We use Proposition \ref{prop_joaquin} to prove preliminary regularity results as usual in the literature. Now, the relevance of estimating $F(X_\varepsilon)-F(Y_\varepsilon)$ in terms of the eigenvalues of $X_\varepsilon-Y_\varepsilon$ becomes clear. Hence, we proceed by verifying that \ref{assump_deg_ellip} is equivalent to supposing that $F$ is degenerate elliptic and Lipschitz continuous. 

\begin{Lemma}\label{lem_murray}
Suppose $F:S(d)\to\mathbb{R}$ satisfies {\rm\ref{assump_deg_ellip}}. Then $F$ is Lipschitz continuous, with constant $C_F$, and monotone non-increasing. That is, for every $M,N \in S(d)$, we have
\begin{equation}\label{eq_F_lip}
	|F(M)- F(N)|\leq C_F\left|N-M\right|,
\end{equation}
and
\begin{equation}\label{eq_F_mono}
	F(M)\leq F(N),
\end{equation}
provided $N\leq M$.
Conversely, suppose $F$ satisfies \eqref{eq_F_lip} and \eqref{eq_F_mono}. Then it also satisfies {\rm\ref{assump_deg_ellip}}.
\end{Lemma}
\begin{proof}
We start by proving that {\rm\ref{assump_deg_ellip}} implies \eqref{eq_F_lip} and \eqref{eq_F_mono}. Indeed, if $M\geq N$ then $(N-M)_+=0$ and we immediately get \eqref{eq_F_mono}. Also, since $\left|(N-M)_+\right|\leq |N-M|$, we get
\begin{align*}
    F(M)-F(N)\leq C_F\left|N-M\right|.
\end{align*}
Swapping $M$ and $N$ we get
\begin{align*}
    F(N)-F(M)\leq C_F\left|N-M\right|,
\end{align*}
and \eqref{eq_F_lip} follows.

Now we prove that  \eqref{eq_F_lip} and \eqref{eq_F_mono} imply {\rm\ref{assump_deg_ellip}}. Fix $M,N\in S(d)$ arbitrarily and recall that $M-N=(M-N)_+-(M-N)_-$. Then 
\begin{align*}
    F(M)=&F(N+(M-N))\leq F(N-(M-N)_-)\leq F(N)+C_F|(N-M)_+|,
\end{align*}
where the first inequality follows from \eqref{eq_F_mono}, and the second one is a consequence of \eqref{eq_F_lip}.
\end{proof}

\section{Interior Lipschitz continuity}\label{sec_lipschitz}

We reduce the problem posed in $\Omega\subset\mathbb{R}^d$ to an equation prescribed in the unit ball, $B_1\subset\mathbb{R}^d$. Let $\Omega'\Subset\Omega$. For every $r\in(0,1)$, one can find a natural number $n=n(r)\in\mathbb{N}$ such that there exists a subset $\{x_1,\ldots,x_n\}\subset \Omega'$, with $\overline{B_r(x_i)}\subset\Omega$ for every $i\in\{1,\ldots,n\}$. In addition, the family $(B_{r/2}(x_i))_{i=1}^n$ covers $\Omega'$; that is,
\[
	\Omega'\subset\bigcup_{i=1}^nB_{r/2}(x_i).
\]
As a result, we suppose $\Omega'=B_{r/2}(x_1)$, with $x_1=0\in\Omega'$, and prescribe our problem of interest in open balls. In what follows, we denote with $C(d)$ any constant depending only on the dimension; this notation refers to possibly different constants within our arguments.

The next lemma accounts for the H\"older regularity of sub-solutions to \eqref{eq_main}. The strategy of the proof yields a modulus of continuity depending explicitly on the growth regime of $H$.

\begin{Lemma}[H\"older continuity for sub-solutions]\label{lem_akiva}
Let $u\in {\rm USC}(B_1)$ be a sub-solution to \eqref{eq_main}. Suppose assumptions {\rm \ref{assump_deg_ellip}}-{\rm\ref{assump_homog_F}} are in force. Suppose further $H$ satisfies {\rm\ref{assump_H_new}} with $m>2$. Then
\begin{align*}
    |u(x)-u(y)|\leq K|x-y|^\gamma,
\end{align*}
for every $x,y\in B_{1/4}$, where
\begin{align}\label{Equation_gamma}
\gamma=\frac{m-2}{m-1},
\end{align}
and
\[
    K:=2^\frac{1}{m-1}\left( 4^\frac{m}{m-1}\left(\frac{C_F C(d)}{C_2\gamma^m}\right)^\frac{1}{m-1}+4\left(\frac{\|f\|_{L^\infty(B_1)}+C_1}{C_2\gamma^m}\right)^\frac{1}{m}\right).
\]

\end{Lemma}
\begin{proof}\label{Lemma_Holder1}
Fix $x\in B_{1/2}$ and define the function $\phi:B_{1/2}(x)\to\mathbb{R}$ as
\[
	\phi(y):=K\left(\frac{1}{4}-|y-x|^2\right)^{-1}|y-x|^\gamma.
\]
To establish the lemma, it suffices to show that 
\begin{equation}\label{eq_omega}
	w(y):=u(y)-u(x)-\phi(y)\leq 0,
\end{equation}
for every $y\in B_{1/2}(x)$. Indeed, it would imply that, for every $y\in B_{1/4}(x)$,
\[
    u(y)-u(x)\leq \phi(y)=K\left(\frac{1}{4}-|y-x|^2\right)^{-1}|y-x|^\gamma\leq \frac{16}{3}K|y-x|^\gamma.
\]
We split the remainder of the proof into three steps.

\medskip

\noindent{\bf Step 1 -} Here, we prove that $w$, as defined in \eqref{eq_omega}, has no local maximum in $B_{1/2}(x)\setminus\{x\}$. Suppose otherwise, and let $y^*\in B_{1/2}(x)\setminus\{x\}$ be a point of local maximum for $w$. Because $u$ is a viscosity sub-solution of \eqref{eq_main}, Definition \ref{israelidates} implies
\[
	F(D^2\phi(y^*))+H(D\phi(y^*),y^*)\leq f(y^*).
\]
We will produce a contradiction by verifying that $\phi$ is also a strict supersolution.

\medskip

\noindent{\bf Step 2 -} Without loss of generality, set $x=0$ and notice that
\[
    D\phi(y)=K\left( \frac{\frac{\gamma}{4} |y|^{\gamma-2}+(2-\gamma)|y|^\gamma}{\left(\frac{1}{4}-|y|^2\right)^2} \right)y,
\]
and
\begin{equation}\label{Estimate_gradiente_Lemma1}
    |D\phi(y)|^m  \geq K^m\left(\frac{1}{4}-|y|^2\right)^{-2m}|y|^{\gamma-2}\frac{\gamma^m}{4^m},
\end{equation}
where the choice of $\gamma$ in \eqref{Equation_gamma} is instrumental. Moreover, a dimensional constant $C(d)>0$ exists such that 
\[
    |D^2\phi(y)|\leq C(d) K\left(\frac{1}{4}-|y|^2\right)^{-3}|y|^{\gamma-2}.
\]
Because $F$ is Lipschitz continuous and satisfies $F(0)=0$, we obtain
\begin{equation}\label{Estimate_Lipschitz_Lemma1}
   - F(D^2\phi(y))\leq C_F |D^2\phi| \leq C_F  C(d) K\left(\frac{1}{4}-|y|^2\right)^{-3}|y|^{\gamma-2}.
\end{equation}
Combining \eqref{Estimate_gradiente_Lemma1} with \eqref{Estimate_Lipschitz_Lemma1}, we get
\[
	\begin{split}
	   	F(D^2\phi(y))+H(D\phi(y),y)&\geq- C_F  C(d) K\left(\frac{1}{4}-|y|^2\right)^{-3}|y|^{\gamma-2}-C_1\\
			&\quad+C_2K^m\left(\frac{1}{4}-|y|^2\right)^{-2m}|y|^{\gamma-2}\gamma^m\\
    			&\geq -C_1+\left(-C_F  C(d) K+ C_2\frac{\gamma^m}{4^m} K^m  \right),
	\end{split}
\]
where the last inequality holds because $(1/4-t^2)^{-3}t^{\gamma-2}>1$, for every $t \in (0,1/2)$ and $\gamma\in (0,1)$. Finally, the choice of $K$ ensures that
\[
    \left(-C_F  C(d) K+ C_2\frac{\gamma^m}{4^m} K^m  \right)> \|f\|_{L^\infty(B_1)}+C_1;
\]
hence, $\phi$ is a supersolution of \eqref{eq_main}, and we obtain a contradiction. Therefore, $w$ does not have an interior local maximum point. In the next step, we prove that $w$ cannot attain its local maximum on $\partial B_{1/2}$.

\medskip

\noindent{\bf Step 3 -} To see that $w$ does not attain a local maximum on $\partial B_{1/2}$, start by noticing that $\phi(y)$ blows up as $y\to\partial B_{1/2}(x)$. As a consequence, the supremum of $w$ in $\overline{B_{1/2}}$ has to be attained in $B_{1/2}$. 

Because of Step 2, it cannot be attained in $B_{1/2}(x)\setminus\{x\}$; hence, the supremum of $w$ is attained at $x$. Because $w(x)=0$, we conclude $w\leq 0$ in $B_{1/2}(x)$, and the proof is complete.
\end{proof}

In the sequel, we produce a Lipschitz-regularity result for the solutions to \eqref{eq_main}. Our argument is based on Proposition \ref{prop_joaquin}, and follows closely the reasoning developed in \cite{Armstrong-Tran_2015}. The main difference stems from the (fully) nonlinear character of the problem. Here, we resort to assumption \ref{assump_deg_ellip} and explore the interplay between the operator $F$, the eigenvalues of a given matrix, and the Hamiltonian $H$. In what follows, we include the sub-quadratic case $1<m\leq 2$.

Because the source term $f$ is Lipschitz-continuous and bounded in $\Omega$, we can absorb it into the Hamiltonian $H$, at the expense of changing the constants appearing in \eqref{yafo}-\eqref{beersheva} accordingly. In doing so, we are allowed to examine the homogeneous variant of \eqref{eq_main} given by
\begin{equation}\label{eq_homogeneous}
	F(D^2u)+H(Du,x)=0\hspace{.2in}\mbox{in}\hspace{.2in}B_1.
\end{equation}
In the sequel we detail the proof of Theorem \ref{teo_main1}.

\begin{proof}[Proof of Theorem \ref{teo_main1}]
We start by setting $L$ as
\begin{equation}\label{eq_benzion}
\begin{aligned}
    L=\,\max\Bigg\{ &2^\frac{1}{m-1}\left( \left(\frac{C_F C(d)}{C_2\gamma^m}\right)^\frac{1}{m-1}+\left(\frac{C_1}{C_2\gamma^m}\right)^\frac{1}{m}\right),\\
  &  \hspace*{2cm}  2\left(3^mC_F C(d)\frac{C_3}{C_2}\right)^\frac{1}{m-1},\,\left(\frac{C_1}{C_3}\right)^\frac{1}{m}\Bigg\}.
\end{aligned}
\end{equation}
It suffices to prove that, for every $\hat{x}\in B_{1/2}$,
\begin{align*}
    \limsup_{x\to \hat{x}}\frac{u(\hat{x})-u(x)}{|\hat{x}-x|}\leq L.
\end{align*}
We suppose there exists $\hat{x}\in B_{1/2}$ such that
\begin{align}\label{Inequality_continuity_contradiction}
    \limsup_{x\to \hat{x}}\frac{u(\hat{x})-u(x)}{|\hat{x}-x|}> L.
\end{align}
By combining Proposition \ref{prop_joaquin} and Lemma \ref{lem_murray} with the conditions in assumptions \ref{assump_deg_ellip}-\ref{assump_H_new}, we obtain a contradiction and complete the proof. For ease of presentation, we split the argument into four steps.

\medskip

\noindent{\bf Step 1 -} Consider first an auxiliary function. Let $\phi:B_{3/4}\to [1,\infty)$ be such that $\phi\equiv 1$ in $B_{1/2}$, with $\phi(x)\to \infty$ as $|x|\to 3/4$. Suppose also 
\begin{equation}\label{Inequality_cutoff_phi}
	|D^j\phi(x)|\leq C(d)\left(\phi(x)\right)^{jm+(1-j)},
\end{equation}
for $j\in\{1,2\}$, $x\in B_{3/4}$, and some dimensional constant $C(d)>0$. 
For $\alpha >0$, denote with $\Psi:B_{3/4}\times B_{3/4}\to\mathbb{R}$ the function
\[
	\Psi(x,y):=u(x)-u(y)-L\phi(y)\left|x-y\right|-\frac{1}{2\alpha}\left|x-y\right|^2.
\]
For $0<\alpha \ll 1$ sufficiently small, we claim that there exist $x_\alpha, y_\alpha \in B_{3/4}$ such that
\begin{equation}\label{eq_designated}
	\Psi(x_\alpha,y_\alpha)=\sup_{x,y\in B_{3/4}}\Psi(x,y)>0.
\end{equation}
 In addition, the function $\phi(y)$ localizes $y_\alpha$ away from the boundary $\partial B_{3/4}$. Also, because $0<\alpha\ll1$, the term $-\frac{1}{2\alpha}\left|x-y\right|^2$ ensures that $x_\alpha$ is close to $y_\alpha$ and, therefore, also away from $\partial B_{3/4}$. Finally, $x_\alpha\neq y_\alpha$, since otherwise the supremum in \eqref{eq_designated} would be zero.

Notice that
\[
    \frac{1}{2\alpha}\left|x_\alpha-y_\alpha\right|^2\leq \osc_{B_{3/4}} u\leq 2\|u\|_{L^\infty(B_{3/4})}.
\]
Hence by the continuity of $u$, we get
\begin{align}\label{eq_sups1}
   \limsup_{\alpha\to 0} &\left( L\phi(y_\alpha)\left|x_\alpha-y_\alpha\right|+\frac{1}{2\alpha}\left|x_\alpha-y_\alpha\right|^2 \right) \nonumber\\
   \leq &\limsup_{\alpha\to 0} \sup \left\{ u(y)-u(z): y,z\in B_{3/4}, |y-z|\leq \left(2\alpha \osc_{B_{3/4}}u\right)^\frac{1}{2} \right\}\\
   =&0\nonumber.
\end{align}
In the case $m>2$, Lemma \ref{lem_akiva} yields
\begin{align*}
    L\phi(y_\alpha)\left|x_\alpha-y_\alpha\right|+\frac{1}{2\alpha}\left|x_\alpha-y_\alpha\right|^2\leq u(x_\alpha)-u(y_\alpha)\leq \tilde{K}\left|x_\alpha-y_\alpha\right|^\gamma,
\end{align*}
where $\gamma=\frac{m-2}{m-1}$ and $\tilde{K}$ stands for the constant $K$ in Lemma \ref{lem_akiva} in the case $f \equiv 0$. Thus,
\begin{align}\label{eq_sups2}
    \phi^{m-1}(y_\alpha)|x_\alpha-y_\alpha|\leq L^{1-m}\tilde{K}^{m-1}.
\end{align}
Next we resort to Proposition \ref{prop_joaquin}. 

\medskip

\noindent{\bf Step 2 -} Because $(x_\alpha,y_\alpha)$ is a maximum point for $\Psi$ and $u$ solves \eqref{eq_homogeneous}, Proposition \ref{prop_joaquin} yields symmetric matrices $X_{\varepsilon,\alpha}$ and $Y_{\varepsilon,\alpha}$ such that
\begin{equation}\label{Inequality_matrix_2}
    \begin{pmatrix}
      X_{\varepsilon,\alpha} & 0\\
      0 & -Y_{\varepsilon,\alpha}
    \end{pmatrix}\leq J_\alpha+\varepsilon J_\alpha^2,
\end{equation}
for every $\varepsilon>0$ and $\alpha>0$, sufficiently small. Moreover, 
\begin{equation}\label{Equation_subjet}
    F(X_{\varepsilon,\alpha})+H\left(P_\alpha,x_\alpha\right)\leq 0\leq F(Y_{\varepsilon,\alpha})+H\left(P_\alpha-Q_\alpha,y_\alpha\right),
\end{equation}
where
\[
    \sigma_\alpha:=\frac{x_\alpha-y_\alpha}{|x_\alpha-y_\alpha|},\hspace{.1in}P_\alpha:=\left(L\phi(y_\alpha)+\frac{|x_\alpha-y_\alpha|}{\alpha}\right)\sigma_\alpha,\hspace{.1in}
\]
and
\[
    Q_\alpha:=L|x_\alpha-y_\alpha|D\phi(y_\alpha).
\]
Finally, we write $J_\alpha$ as
\[
    J_\alpha=\frac{L\phi(y_\alpha)}{|x_\alpha-y_\alpha|}\begin{pmatrix}
      Z_1 & -Z_1\\
      -Z_1 & Z_1
    \end{pmatrix} +\frac{1}{\alpha}\begin{pmatrix}
      I & -I\\
      -I & I
    \end{pmatrix}+L\begin{pmatrix}
      0 & Z_2 \\
      Z_2^T & Z_3
    \end{pmatrix},
\]
with $Z_1:=I-\sigma_\alpha\otimes \sigma_\alpha$, $Z_2:=D\phi(y_\alpha)\otimes \sigma_\alpha$, and 
\[
	Z_3:=-(Z_2+Z_2^T)+D^2\phi(y_\alpha)|x_\alpha-y_\alpha|.
\]
In the next step, we estimate $F(Y_{\varepsilon,\alpha})-sF(X_{\varepsilon,\alpha})$ from above.

\medskip

\noindent{\bf Step 3 -} It follows from {\rm\ref{assump_deg_ellip}} that
\begin{equation}\label{eq_murray}
    F(Y_{\varepsilon,\alpha})-F(sX_{\varepsilon,\alpha})\leq C_F\left|\left(sX_{\varepsilon,\alpha}-Y_{\varepsilon,\alpha}\right)_+\right|.
\end{equation}
For $s>0$, let $A_s$ be given by
\[
    A_s=\begin{pmatrix}
      s^2 I & sI\\
      sI & I
    \end{pmatrix}.
\]
Multiply both sides of \eqref{Inequality_matrix_2} by $A_s$ and evaluate the resulting inequality at vectors of the form $(\omega, \omega)\in \mathbb{R}^{2d}$. As a consequence, one obtains
\[
    \omega^T\left((s^2+s)X_{\varepsilon,\alpha}-(s+1)Y_{\varepsilon,\alpha}\right)\omega\leq L(s+1)\omega^T\left(Z_2+sZ_2^T+Z_3\right)\omega+O(\varepsilon).
\]
Set $s:=1+\beta|x_\alpha-y_\alpha|$, with $\beta=\overline{\beta}\phi^{m-1}(y_\alpha)$, for $\overline{\beta}$ yet to be fixed. It follows that 
\begin{align*}
   \omega^T\left(sX_{\varepsilon,\alpha}-Y_{\varepsilon,\alpha}\right)\omega&\leq L\omega^T\left(Z_2+sZ_2^T+Z_3\right)\omega+O(\varepsilon)\\
    &= L\omega^T\left((s-1)\sigma_\alpha \otimes D\phi(y_\alpha)\right)\omega\\
	&\quad+L\omega^T\left(D^2\phi(y_\alpha)|x_\alpha-y_\alpha|\right)\omega+O(\varepsilon)\\
    &\leq|\omega|^2L\left((s-1)|D\phi(y_\alpha)||x_\alpha-y_\alpha|\right)\\
	&\quad+|\omega|^2L\left(|D^2\phi(y_\alpha)||x_\alpha-y_\alpha|\right)+O(\varepsilon)\\
    &\leq|\omega|^2L C(d)\overline{\beta}\phi^{2m-1}(y_\alpha)|x_\alpha-y_\alpha| +O(\varepsilon).
\end{align*}
In conclusion, 
\begin{align}\label{Inequality_RHS}\nonumber
    F(Y_{\varepsilon,\alpha})-sF(X_{\varepsilon,\alpha})&\leq C_F\left|\left(sX_{\varepsilon,\alpha}-Y_{\varepsilon,\alpha}\right)_+\right|\\
    	&\leq C_FC(d)L\overline{\beta}\phi^{2m-1}(y_\alpha)|x_\alpha-y_\alpha|+O(\varepsilon).
\end{align}
In what follows, we estimate $F(Y_{\varepsilon,\alpha})-sF(X_{\varepsilon,\alpha})$ from below.

\medskip

\noindent{\bf Step 4 -} We start with three auxiliary inequalities. Because of \eqref{Inequality_cutoff_phi}, we have
\begin{equation}\label{Inequality_Qalpha}
    \frac{|Q_\alpha|}{|x_\alpha-y_\alpha|}\leq C(d)\phi^{m-1}(y_\alpha)\leq C(d)L^{1-m}|P_\alpha|^m.
\end{equation}
More generally, for any $\theta>0$,
\begin{equation}\label{Inequality_Qalpha2}
    \frac{|Q_\alpha|}{|x_\alpha-y_\alpha|}\leq C(d) L^{1-\theta} \phi^{m-\theta}(y_\alpha)|P_\alpha|^\theta.
\end{equation}
Also,
\begin{equation}\label{Inequality_Palpha}
    L|D^2\phi(y_\alpha)|\leq C(d) L^{1-m}\phi^{m-1}(y_\alpha)|P_\alpha|^m.
\end{equation}

In the sub-quadratic case $1<m\leq 2$, one combines \eqref{Inequality_Qalpha2} with $\theta=1$ and \eqref{eq_sups1} to get
\begin{align*}
    \lim_{\alpha\to 0} \frac{|Q_\alpha|}{|P_\alpha|}\leq \lim_{\alpha\to 0} C(d) |\phi(y_\alpha)|^{m-1}|x_\alpha-y_\alpha|=0.
\end{align*}
For $m>2$, in the superquadratic case, \eqref{Inequality_Qalpha} builds upon \eqref{eq_sups2} to produce
\begin{align*}
    |Q_\alpha|\leq C(d) L^{1-m}\tilde{K}^{m-1}|P_\alpha|.
\end{align*}
Hence, in either case, we have $|Q_\alpha|\leq |P_\alpha|$, for $\alpha$ small enough, since \eqref{eq_benzion} implies
\begin{align}\label{Inequality_lowerbound_L}
    L\geq 2^\frac{1}{m-1}\left( \left(\frac{C_F C(d)}{C_2\gamma^m}\right)^\frac{1}{m-1}+\left(\frac{C_1}{C_2\gamma^m}\right)^\frac{1}{m}\right).
\end{align}

Now we use \eqref{Equation_subjet} to write
\[
	\begin{split}
    F(Y_{\varepsilon,\alpha})-F(sX_{\varepsilon,\alpha})&\geq sH(P_\alpha,x_\alpha)-H(P_\alpha-Q_\alpha,y_\alpha)\\
    	&\geq(s-1)H(P_\alpha,x_\alpha)-C_3(1+2|P_\alpha|)^{m-1}|Q_\alpha|\\
	&\quad-(C_3|P_\alpha|^m+C_1)|x_\alpha-y_\alpha|\\
           &\geq(s-1)(C_2|P_\alpha|^m-C_1)-C_3(1+2|P_\alpha|)^{m-1}|Q_\alpha|\\
	&\quad-(C_3|P_\alpha|^m+C_1)|x_\alpha-y_\alpha|.
\end{split}
\]
We used \ref{assump_homog_F} in the first inequality, whereas $|Q_\alpha|\leq |P_\alpha|$ leads to the second one. 
Since $|P_\alpha|\geq L>1$, it follows that  $1+2|P_\alpha|\leq 3|P_\alpha|$, and from the lower bound $L\geq (C_1/C_3)^{1/m}$, we also obtain $C_3|P_\alpha|^m+C_1\leq 2C_3|P_\alpha|^m.$
Thus, we can further estimate
\[
	\begin{split}
		& F(Y_{\varepsilon,\alpha})-F(sX_{\varepsilon,\alpha})\\
			&\geq (s-1)C_2|P_\alpha|^m-C_3(3|P_\alpha|)^{m-1}|Q_\alpha|-2C_3|P_\alpha|^m|x_\alpha-y_\alpha|\\
			&\geq\left(\beta C_2|P_\alpha|^m-C_3(3|P_\alpha|)^{m-1}\frac{|Q_\alpha|}{|x_\alpha-y_\alpha|}\right)|x_\alpha-y_\alpha|\\
			&\quad-2C_3|P_\alpha|^m|x_\alpha-y_\alpha|\\
		    	&\geq \left(\beta C_2|P_\alpha|^m-C_3(3|P_\alpha|)^{m-1}\left(\phi^{m-1}(y_\alpha)|P_\alpha| \right)\right)|x_\alpha-y_\alpha|\\
			&\quad-2C_3|P_\alpha|^m|x_\alpha-y_\alpha|\\	
			&=|P_\alpha|^m\left(\beta C_2-C_33^{m-1}\phi^{m-1}(y_\alpha) -2C_3\right)|x_\alpha-y_\alpha|.
	\end{split}
\]
We used \eqref{Inequality_Qalpha2} with $\theta=1$ in the last inequality. Because $\beta=\overline{\beta}\phi^{m-1}(y_\alpha)$, we get
\begin{align}\label{Inequality_LHS}
    F(Y_{\varepsilon,\alpha})-F(sX_{\varepsilon,\alpha})\geq& |P_\alpha|^m\phi^{m-1}\left(\overline{\beta} C_2-C_33^{m-1} -2C_3\phi^{1-m}(y_\alpha)\right)|x_\alpha-y_\alpha|\nonumber\\
    \geq &|P_\alpha|^m\phi^{m-1}|x_\alpha-y_\alpha|,
\end{align}
provided we set $\overline{\beta}=\frac{C_3}{C_2}(3^{m-1}+2)$.

Combining \eqref{Inequality_LHS} and \eqref{Inequality_RHS} we get
\[
    |P_\alpha|^m\phi^{m-1}|x_\alpha-y_\alpha|\leq C_F\left( C(d)L\overline{\beta}\phi^{2m-1}(y_\alpha)|x_\alpha-y_\alpha|+O(\varepsilon)\right).
\]
Let $\varepsilon\to 0$ and divide both sides of the former inequality by the quantity $|x_\alpha-y_\alpha|\phi^{2m-1}(y_\alpha)$. Then
\[
    L^m\leq \frac{|P_\alpha|^m}{\phi^m(y_\alpha)}\leq \overline{C} L,
\]
where $\overline{C}=C_F C(d) \overline{\beta}$. This is a contradiction since
\[
    L\geq 2\left(3^mC_F C(d)\frac{C_3}{C_2}\right)^\frac{1}{m-1}.
\]
Therefore we have proven that $u$ is locally Lipschitz continuous, with constant $L$ given by \eqref{eq_benzion}.
\end{proof}

\begin{Remark}
The proof of Theorem \ref{teo_main1} provides a constructive way to produce the Lipschitz constant $C$ associated with $u$. In fact, this is given by $L>0$, as defined in \eqref{eq_benzion}.
\end{Remark}

\section{A two-phase free boundary problem}\label{cleese}

Now, we explore a consequence of Lemma \ref{lem_akiva} in the context of free boundary problems. It concerns the existence of a viscosity solution to 
\begin{equation}\label{aaronjames}
	\begin{cases}
		-{\rm Tr}(A(x)D^2u)+H(Du,x)=\lambda_+\chi_{\{u>0\}}+\lambda_-\chi_{\{u<0\}}&\hspace{.2in}\mbox{in}\hspace{.2in}\Omega(u)\\
		u=g&\hspace{.2in}\mbox{on}\hspace{.2in}\partial\Omega
	\end{cases}
\end{equation}
where $\Omega$ is a $C^2$-domain, $A:\Omega \to S(d)$ is degenerate elliptic, $0<\lambda_-<\lambda_+$ are constants, $g\in C^{0,\frac{m-2}{m-1}}(\partial\Omega)$ is given, and $\Omega(u)$ is given by
\[
	\Omega(u):=\left\lbrace x\in\Omega \,|\, u(x)\neq 0\right\rbrace.
\]
We notice the equation holds \emph{only} where the solution does not vanish, and hence, no information is available across the free boundary $\Gamma(u):=\partial\{u>0\}\cup\partial\{u<0\}$.

We prove the existence of a locally H\"older-continuous viscosity solution to \eqref{aaronjames} with suitable, estimates. To do that, we introduce an assumption on the data $A:\Omega\to S(d)$ and $g:\partial\Omega\to\mathbb{R}$.	

\begin{assump}\label{assump_g}
We suppose $A:\Omega\to S(d)$ to be degenerate elliptic and bounded from above. In addition, there exists $\lambda>0$ such that 
\[
	\nu(x)^TA(x)\nu(x)\geq\lambda
\]
for every $x\in\partial\Omega$. Also, we suppose $g\in C^{0,\frac{m-2}{m-1}}(\partial\Omega)$, with
\[
	|g(x)-g(y)|\leq K|x-y|^\frac{m-2}{m-1}
\]
for every $x,y\in\partial\Omega$, where $K>0$ is fixed, though yet to be determined. In addition, suppose
\[
	0<\inf_{x\in \partial\Omega}g(x)< 2\lambda_-.
\]
\end{assump}

The importance of \ref{assump_g} is in unlocking an intermediate step in our analysis, namely \cite[Theorem 2.12]{Dolcetta-Leoni-Porretta_2010}. In fact, the superquadratic character of \eqref{aaronjames} introduces a number of subtleties in the arguments leading to the existence of solutions. See the discussion in \cite[Section 2.3]{Dolcetta-Leoni-Porretta_2010}.

\begin{Theorem}[Existence of solutions]\label{paulajones} Let $\Omega\subset\mathbb{R}^d$ be an open, bounded domain of class $C^2$. Suppose assumption {\rm\ref{assump_g}} is in force. Then there exists a viscosity solution $u\in C(\Omega)$ to the problem \eqref{aaronjames}. In addition, we have $u\in C^{0,\frac{m-2}{m-1}}_{\rm loc}(\Omega)$. Finally, for every $\Omega'\Subset\Omega$, there exists a positive constant $C=C \left( m, \|A\|_{L^\infty(\Omega)}, K, \lambda, {\rm diam}(\Omega),{\rm dist}(\Omega',\partial\Omega) \right)$ such that 
\[
	\|u\|_{C^{0,\frac{m-2}{m-1}}(\Omega')}\leq C\left(1+\|u\|_{L^\infty(\Omega)}+{\rm max}\{|\lambda_+|,|\lambda_-|\}\right).
\]
\end{Theorem}

The proof of Theorem \ref{paulajones} combines several ingredients. First, we consider a family of auxiliary equations indexed by a parameter $\varepsilon>0$. For each equation in the family, the existence of a (unique) viscosity solution follows from \cite[Theorem 2.12]{Dolcetta-Leoni-Porretta_2010}. Lemma \ref{lem_akiva} implies estimates independent of $\varepsilon>0$ and allows us to apply Schauder's Fixed Point Theorem to conclude the argument. We proceed by introducing an auxiliary problem.

For $v\in C(\overline{\Omega})$ and $0<\varepsilon<1$, define $g_\varepsilon^v:\mathbb{R}^d\to\mathbb{R}$ as
\[
	g_\varepsilon^v(x):=\max\left(\min\left(\frac{v(x)+\varepsilon}{2\varepsilon},1\right),0\right),
\]
if $x\in\Omega$, and $g_\varepsilon^v\equiv 0$ in $\mathbb{R}^d\setminus\Omega$. Now, for $x\in\Omega$, let $h_\varepsilon^v(x):=\left(g_\varepsilon^v\ast\eta_\varepsilon\right)(x)$, where $\eta_\varepsilon$ is a standard mollifier. We consider the auxiliary problem
\begin{equation}\label{tripp}
	\begin{cases}
		\varepsilon u-{\rm Tr}(A(x)D^2u)+H(Du,x)=\lambda_+h_\varepsilon^v+\lambda_-(1-h_\varepsilon^v)&\hspace{.2in}\mbox{in}\hspace{.2in}\Omega\subset\mathbb{R}^d\\
		u=g&\hspace{.2in}\mbox{on}\hspace{.2in}\partial\Omega.
	\end{cases}
\end{equation}
Through a fixed-point approach, we show the existence of a solution $u_\varepsilon\in C(\Omega)$ to the Dirichlet problem
\begin{equation}\label{kenstarr}
	\begin{cases}
		\varepsilon u_\varepsilon-{\rm Tr}(A(x)D^2u_\varepsilon)+H(Du_\varepsilon,x)=\lambda_+h_\varepsilon^{u_\varepsilon}+\lambda_-(1-h_\varepsilon^{u_\varepsilon})&\hspace{.2in}\mbox{in}\hspace{.2in}\Omega\subset\mathbb{R}^d\\
		u_\varepsilon=g&\hspace{.2in}\mbox{on}\hspace{.2in}\partial\Omega.
	\end{cases}
\end{equation}
This argument relies on a variant of Lemma \ref{lem_akiva} applied to \eqref{tripp}. By taking the limit $\varepsilon\to 0$ in \eqref{kenstarr} and applying Lemma \ref{lem_akiva} once more, we obtain the existence of a viscosity solution to \eqref{aaronjames}. The first step towards proving Theorem \ref{paulajones} is the following proposition.

\begin{Proposition}\label{lemur}
Let $\Omega\subset\mathbb{R}^d$ be an open, bounded domain of class $C^2$. Suppose assumption {\rm\ref{assump_g}} is in force. Then, for every $0<\varepsilon<1/4$, there exists $u_\varepsilon\in C(\Omega)$ solving \eqref{kenstarr} in the viscosity sense. In addition, $u_\varepsilon\in C^{0,\frac{m-2}{m-1}}_{\rm loc}(\Omega)$. Moreover, for every $\Omega'\Subset\Omega$, there exists a positive constant $C=C \left( m, \|A\|_{L^\infty(\Omega)}, K, \lambda, {\rm diam}(\Omega),{\rm dist}(\Omega',\partial\Omega) \right)$ such that
\begin{equation}\label{eilat}
	\|u\|_{C^{0,\frac{m-2}{m-1}}(\Omega')}\leq C\left(1+\|u\|_{L^\infty(\Omega)}+{\rm max}\{|\lambda_+|,|\lambda_-|\}\right).
\end{equation}
\end{Proposition}
\begin{proof}
For ease of presentation, we split the proof into three steps.

\medskip

\noindent{\bf Step 1 - } Notice that, given $v\in C(\Omega)$, we have
\[
	|\lambda_+h_\varepsilon^v(x)+\lambda_-(1-h_\varepsilon^v(x))-\lambda_+h_\varepsilon^v(y)-\lambda_-(1-h_\varepsilon^v(y))|\leq \omega_{v,\varepsilon}(|x-y|),
\]
where $\omega_{v,\varepsilon}(\cdot)$ is a modulus of continuity depending on $v$ and $\varepsilon>0$. Hence, the right-hand side of the equation in \eqref{tripp} is a continuous function up to the boundary $\partial\Omega$. A straightforward application of \cite[Theorem 2.12]{Dolcetta-Leoni-Porretta_2010} ensures the existence of a unique viscosity solution $u^v_\varepsilon$ to \eqref{tripp}.

Notice also the proof of Lemma \ref{lem_akiva} extends to the case of \eqref{tripp}. As a consequence, $u_\varepsilon^v\in C^{0,\frac{m-2}{m-1}}_{\rm loc}(\Omega)$ for every $0<\varepsilon<1/4$ and every $v\in C(\overline{\Omega})$. Moreover, for every $\Omega'\Subset\Omega$, there exists a constant $C>0$ depending on the data of the problem and $\Omega'$, but not depending on $\varepsilon$ or $v$, such that
\begin{equation}\label{hillary}
	\left\|u_\varepsilon^v\right\|_{C^{0,\frac{m-2}{m-1}}(\Omega')}\leq C.
\end{equation}
uniformly in $v\in C(\Omega)$ and $\varepsilon\in(0,1/4)$.

\medskip

\noindent{\bf Step 2 - }We now define $\mathcal{K}\subset C(\Omega)$ as
\[
	\mathcal{K}:=\left\lbrace w\in C(\overline{\Omega})\; : \; \|w\|_{L^\infty (\Omega)} \leq C_0\right\rbrace,
\]
where $C_0>0$ will be chosen later. Notice that $\mathcal{K}$ is closed in $C(\overline{\Omega})$. In the sequel, we define a map $T:\mathcal{K}\to C(\overline{\Omega})$. For fixed $\varepsilon\in(0,1/4)$, take $v\in \mathcal{K}$ and denote with $u_\varepsilon^v$ the unique solution to \eqref{tripp}, whose existence follows from the previous step. Define $Tv:=u_\varepsilon^v$ and notice that the existence of a fixed point for $T$ is tantamount to the existence of solutions to \eqref{kenstarr}.

To prove the existence of a fixed point for $T$, we start by noticing that it is possible to choose $C_0>0$, independent on $v$ and $\varepsilon$, such that $T(\mathcal{K})\subset \mathcal{K}$. This follows from the construction of sub and supersolutions for the problem; see the proof of \cite[Theorem 2.12]{Dolcetta-Leoni-Porretta_2010}.

We continue by proving that $T(\mathcal{K})$ is pre-compact. Let $(Tv_n)_{n\in\mathbb{N}}$ be a sequence of elements in $T(\mathcal{K})$. Because $\left| Tv_n \right| \leq C_0$, for every $n\in\mathbb{N}$, the sequence is equibounded. In addition, \eqref{hillary} ensures it is also equicontinuous. Hence, $(Tv_n)_{n\in\mathbb{N}}$ converges to an element $v^*\in T(\mathcal{K})$, through a subsequence if necessary.

Finally, we verify that $T$ is sequentially continuous. Suppose $(v_n)_{n\in\mathbb{N}}\subset \mathcal{K}$ converges to some $v\in\mathcal{K}$. We prove that $Tv_n\to Tv$, as $n\to\infty$. Indeed, because $T(\mathcal{K})$ is pre-compact, we infer that 
$(Tv_n)_{n\in\mathbb{N}}$ converges, through a subsequence if necessary, to some $w\in\mathcal{K}$. The stability of viscosity solutions and the uniqueness available for \eqref{tripp} ensure that $Tv=w$. To conclude that $T$ is continuous, we must verify the former equality does not depend on the particular subsequence. Indeed, suppose a different subsequence yields $Tv_{n_k}\to w'$, as $k\to\infty$. Since $v_{n_k}\to v$ as $k\to\infty$, we reason as before (resorting to the stability of viscosity solutions and the uniqueness for \eqref{tripp}) to obtain $Tv=w'$.

\medskip

\noindent{\bf Step 3 - } The properties of the subset $\mathcal{K}$ and the operator $T$ allow us to apply the Schauder Fixed Point Theorem, as in \cite[Corollary 11.2]{Gilbarg-Trudinger_2001}, to conclude the existence of $u_\varepsilon\in\mathcal{K}$ such that $Tu_\varepsilon=u_\varepsilon$. That is, $u_\varepsilon$ solves \eqref{kenstarr}. Because the conclusions of Lemma \ref{lem_akiva} apply, we have $u_\varepsilon\in C^{0,\frac{m-2}{m-1}}_{\rm loc}(\Omega)$ and the estimate in \eqref{eilat} holds.
\end{proof}

Now, we detail the proof of Theorem \ref{paulajones}.

\begin{proof}[Proof of Theorem \ref{paulajones}]
We take a sequence $(\varepsilon_n)_{n\in\mathbb{N}}$ such that $\varepsilon_n\to 0$ as $n\to\infty$ and consider the sequence $(u_n)_{n\in\mathbb{N}}$ of solutions to
\[
	\varepsilon_nu_n-{\rm Tr}(A(x)D^2u_n)+H(Du_n,x)=\lambda_+h_{\varepsilon_n}^{u_n}+\lambda_-(1-h_{\varepsilon_n}^{u_n})\hspace{.2in}\mbox{in}\hspace{.2in}\Omega.
\]
Because the estimate in \eqref{eilat} holds for $u_n$, for every $n\in\mathbb{N}$, we conclude there exists $u^*\in C^\beta_{\rm loc}(\Omega)$ such that $u_n\to u^*$ in the $C^\beta$-topology, for every $0<\beta<\frac{m-2}{m-1}$. 

Now, let $x\in\{u^*>0\}$ and write $\tau:=u^*(x)$. Suppose $u^*-\varphi$ has a (strict) local maximum at $x$, for $\varphi\in C^2(\Omega)$. There exists a sequence $(x_n)_{n\in\mathbb{N}}$ such that $x_n\to x$ and $u_n-\varphi$ has a local maximum at $x_n$. On the other hand, there exists $N\in\mathbb{N}$ such that 
\[
	u_n(x_n)>\frac{\tau}{2}>\varepsilon_n
\]
provided $n>N$. Hence, 
\[
	\varepsilon_nu_n(x_n)-{\rm Tr}(A(x_n)D^2\varphi(x_n))+H(D\varphi(x_n),x_n)\leq\lambda_+h_{\varepsilon_n}^{u_n}(x_n)-\lambda_-(1-h_{\varepsilon_n}^{u_n}(x_n));
\]
by taking the limit $n\to\infty$, we obtain
\[
	-{\rm Tr}(A(x)D^2\varphi(x))+H(D\varphi(x),x)\leq\lambda_+.
\]
Conversely, suppose $x\in\{u^*<0\}$ and write $\sigma:=u^*(x)$. If $u^*-\varphi$ has a (strict) local maximum at $x$, for $\varphi\in C^2(\Omega)$ we reason as before to conclude
\[
	-{\rm Tr}(A(x)D^2\varphi(x))+H(D\varphi(x),x)\leq-\lambda_-.
\] 
It ensures that $u^*$ is a sub-solution to \eqref{aaronjames} in $\Omega(u)$. An analogous argument ensures that $u^*$ is also a supersolution and completes the proof.
\end{proof}

\begin{Remark}
Our proof of Theorem \ref{paulajones} yields further information since it produces two viscosity inequalities satisfied by the solution \emph{in the region $\{u=0\}$}. Indeed, the viscosity solution to \eqref{aaronjames}, whose existence follows from Theorem \ref{paulajones}, satisfies 
\[
	-|\lambda_-|\leq - {\rm Tr}(A(x)D^2u)+H(Du,x)\leq |\lambda_+|\hspace{.2in}\mbox{in}\hspace{.2in}\Omega
\]
in the viscosity sense. Besides solving the equation in the positive and negative phases, it also solves a pair of inequalities in the whole domain.
\end{Remark}

\bigskip

{\small \noindent{\bf Acknowledgments.} DJ was supported by Funda\c c\~ao para a Ci\^encia e a Tecnologia, through scholarship PD/BD/150354/2019, under POCH funds, co-financed by the European Social Fund and Portuguese National Funds from MCTES, and by the Centre for Mathematics of the University of Coimbra (UIDB/00324/2020, funded by the Portuguese Government through FCT/MCTES). 

\noindent EP partially supported by the Centre for Mathematics of the University of Coimbra (UIDB/00324/2020, funded by the Portuguese Government through FCT/MCTES) and by FAPERJ (grants E26/200.002/2018 and E26/201.390/2021). 

\noindent JMU partially supported by the King Abdullah University of Science and Technology (KAUST), by Funda\c c\~ao para a Ci\^encia e a Tecnologia, through project PTDC/MAT-PUR/28686/2017, and by the Centre for Mathematics of the University of Coimbra (UIDB/00324/2020, funded by the Portuguese Government through FCT/MCTES).}

\end{document}